Australian Journal of Basic and Applied Sciences, 10(9) May 2016, Pages: 8-18

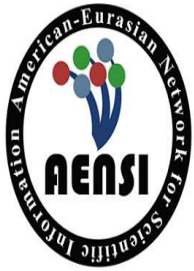

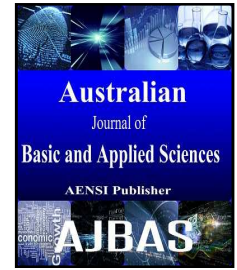

# General Asymptotic Regional Gradient Observer

[1]Asmaa N. Al-Janabi, [2]Raheam A. Al-Saphory, [3]Naseif J. Al-Jawari

[1]Department of Mathematics, College of Science, Al-Mustansriyah University, Baghdad, IRAQ.
[2]Department of Mathematics, College of Education for Pure Sciences, Tikrit University, Tikrit, IRAQ.
[3]Department of Mathematics, College of Science, Al-Mustansriyah University, Baghdad, IRAQ.

**Address For Correspondence:**
Asmaa N. Al-Janabi, 1Department of Mathematics, College of Science, Al-Mustansriyah University, Baghdad, IRAQ.
E-mail: asm_no_2006@yahoo.com

**A B S T R A C T**
The main purpose of this paper is to study and characterize the existing of general asymptotic regional gradient observer which observe the current gradient state of the original system in connection with gradient strategic sensors. Thus, we give an approach based to Luenberger observer theory of linear distributed parameter systems which is enabled to determinate asymptotically regional gradient estimator of current gradient system state. More precisely, under which condition the notion of asymptotic regional gradient observability can be achieved. Furthermore, we show that the measurement structures allows the existence of general asymptotic regional gradient observer and we give a sufficient condition for such asymptotic regional gradient observer in general case. We also show that, there exists a dynamical system for the considered system is not general asymptotic gradient observer in the usual sense, but it may be general asymptotic regional gradient observer. Then, for this purpose we present various results related to different types of sensor structures, domains and boundary conditions in two dimensional distributed diffusion systems.

## INTRODUCTION

The notion of observer theory was introduced for finite dimensional linear system by Luenberger and has been extended to multi-variables systems in (Luenberger, 1966). This theory is generalized to infinite dimensional linear systems characterized by strongly continuous semi-group operatorsby(Gressang and Lamont, 1975). The study of observer concept via the sensors and actuators structure was developed by El Jai and Pritchard in (El Jai and Pritchard 1988). The purpose of an observer is to provide an state estimation for the considered system state as in (El Jai andAmouroux1988). The concept of asymptotic regional construction has been introduced by Al-Saphory and El Jai in (Al-Saphory, 2002), (Al-Saphory and El Jai, 2002), (Al-Saphory, 2011) consists in studying the asymptotic behavior of the systems not in the whole the domain but in an internal sub-region $\omega$ of a spatial domain $\Omega$. The main reason for introducing the concept is motivated by certain concrete-real problems, in mechanic, thermic, environment in(El Jai, *et al.*,1995),(Burns *et al.*, 2009)and (Burns *et al.*, 2010). Commercial buildings are responsible for a significant fraction of the energy consumption and greenhouse gas emissions in the U.S. and worldwide. Consequently, the design, optimization and control of energy efficient buildings can have a tremendous impact on energy cost and greenhouse gas emission. Mathematically, building models are complex, multi-scale, multi-physics, highly uncertain dynamical systems with wide varieties of disturbances(Burns *et al.*, 2009).

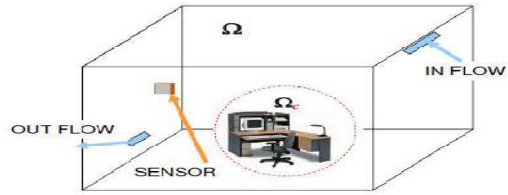

**Fig. 1:** Room control model with sensor, in flow and out flow

In this paper we use a model problem to illustrate that distributed parameter control based on PDEs, combined with high performance computing can be used to provide practical insight into important issues such as optimal sensor/actuator placement (may be best or strategic sensors/ actuators) and optimal supervisory building control. In order to illustrate some of the ideas, we consider the problem illustrated by a single room shown in (Figure 1). This model one can reformulated (Burns *et al*., 2010) as spatial case of more general model of distributed parameter systems and represented in the next section (see Figure 2).

In this paper, we introduced an approach which allow to construct a general asymptotic $\omega_G$-observers in a given sub-region $\omega$ of the domain $\Omega$, in connection with sensors structures this approach based on the concept of asymptotic regional gradient detectability as in (Al-Saphory and El Jai, 2002) and (Al-Saphory *et al*., 2010).

*This paper is organized as follows:*

Section 2 is devoted to the problem statement and preliminaries. We recall that the definitions of $\omega_G$- stability and asymptotic $\omega_G$- detectability and we give some definitions about asymptotic $\omega_G$- observers. Section 3, we have given a characterize for existing a general asymptotic regional gradient observer to provide a general asymptotic regional gradient estimator of gradient state for an original system via sensors structures and we show that there exist a general asymptotic $\omega_G$-observers in the subregion is not general asymptotic $G$- observer in the whole the domain $\Omega$.

*2. Problem Statement:*

Let $\Omega$ be a regular bounded open subset of $R^n$, with smooth boundary $\partial\Omega$ and $\omega$ be subregion of $\Omega$, [0,T], $T > 0$ be a time measurement interval. we denoted $Q = \Omega \times ]0,T[$, $\Sigma = \partial\Omega \times ]0,T[$. We considered distributed parabolic system is described by the following equation:

$$\begin{cases} \frac{\partial x}{\partial t}(\xi,t) = Ax(\xi,t) + Bu(t) & Q \\ x(\xi,0) = x_0(\xi) & \Omega \\ x(\eta,t) = 0 & \Sigma \end{cases} \quad (1)$$

Augmented with the output function
$$y(.,t) = Cx(,t) \quad (2)$$

We have
$$A = \sum_{i,j=1}^{n} \frac{\partial}{\partial x_j}\left(a_{ij}\frac{\partial}{\partial x_j}\right)$$

With $a_{ij} \in \mathcal{D}(\bar{Q})$. Suppose that $-A$ is elliptic, i.e., there exits $\alpha > 0$ such that
$\sum_{i,j=1}^{n} a_{ij}\xi_i\xi_j \geq \alpha \sum_{j=1}^{n}|\xi_j|^2$ a.e. on $Q$, $\forall \xi = (\xi_1,\dots,\xi_n) \in R^n$

Where $A$ is a second order linear differential operator, which generator a strongly continuous semi-group $(S_A(t))_{t\geq 0}$ on the Hilbert space $X$ and is self-adjoin with compact resolvent. The operator $B \in L(R^p, X)$ and $C \in L(R^q, X)$, depend on the structure of actuators and sensors (El Jai and Pritchard 1988). The space $X, U$ and $\mathcal{O}$ be separable Hilbert spaces where $X$ is the state space, $U = L^2(0,T,R^p)$ is the control space and $\mathcal{O} = L^2(0,T,R^q)$ is the observation space where $p$ and $q$ are the numbers of actuators and sensors (see Figure 2).

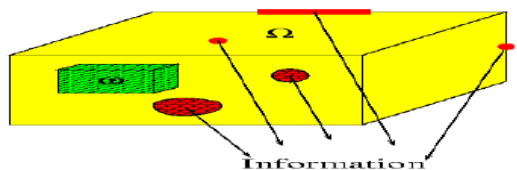

**Fig. 2:** The domain of $\Omega$, the sub-region $\omega$, various sensors locations

Under the given assumption, the system (1) has a unique solution (Luenberger, 1966):
$$x(\xi,t) = S_A(t)x_0(\xi) + \int_0^t S_A(t-\tau)Bu(\tau)d\tau \tag{3}$$

The measurements are obtained through the output function by using of zone, point wise which may located in $\Omega$ (or $\partial\Omega$). (El Jai and Pritchard 1988)
$$y(.,t) = Cx(\xi,t) \tag{4}$$

- We first recall a sensors is defined by any couple $(D, f)$, where $D$ is its spatial support represented by a nonempty part of $\bar{\Omega}$ and $f$ represents the distribution of the sensing measurements on $D$. Depending on the nature of $D$ and $f$, we could have various type of sensors. A sensor may be pointwise if $D=\{b\}$ with $b \in \bar{\Omega}$ and $f = \delta(.-b)$, where $\delta$ is the Dirac mass concentrated at $b$. In this case the operator $C$ is unbounded and the output function (2) can be written in the form
$$y(t) = x(b,t)$$

It may be zonal when $D \subset \bar{\Omega}$ and $f \in L^2(D)$. The output function (2) can be written in the form
$$y(t) = \int_D x(\xi,t)f(\xi)d\xi$$

In the case of boundary zone sensor, we consider $D_i = \Gamma_i \subset \partial\Omega$ and $f_i \in L^2(\Gamma_i)$, the output function (2) can be written as
$$y(.,t) = Cx(.,t) = \int_{\Gamma_i} x(\eta,t)f_i(\eta)d\eta$$

- We define the operator
$$K: x \in X \to Kx = CS_A(.)x \in \mathcal{O}$$

We note that $K^*: \mathcal{O} \to X$ is the adjoint operator of $K$ defined by
$$K^*y^* = \int_0^t S_A^*(s)C^*y^*(s)ds$$

- Consider the operator
$$\nabla: H^1(\Omega) \to (H^1(\Omega))^n$$
$$x \to \nabla x = (\frac{\partial x}{\partial \xi_1}, \dots, \frac{\partial x}{\partial \xi_n})$$

and the adjoint denotes by $\nabla^*$ is given as
$$\nabla^*: (H^1(\Omega))^n \to H^1(\Omega)$$
$$x \to \nabla^* x = v$$

where $v$ is a solution of the Dirichlet problem
$$\begin{cases} \Delta v = -\text{div}(x) & \text{in } \Omega \\ v = 0 & \text{on } \partial\Omega \end{cases}$$

- For a nonempty subset $\omega \subset \Omega$ consider the regional gradient restriction operator
$$\chi_\omega: (H^1(\Omega))^n \to (H^1(\omega))^n$$
$$x \to \chi_\omega x = x|_\omega$$

It's adjoint is denoted by $\chi_\omega^*$.

- $\tilde{\chi}_\omega: \begin{cases} H^1(\Omega) \to H^1(\omega) \\ x \to \tilde{\chi}_\omega x = x|_\omega \end{cases}$

where $x|_\omega$ is the restriction of $x$ to $\omega$.
It's adjoint is denoted by $\tilde{\chi}_\omega^*$.

- Finally, we introduced the operator $H = \chi_\omega \nabla K^*$ from $\mathcal{O}$ into $(H^1(\omega))^n$.

The problem is how to build an approach which observe (estimates) an asymptotic regional gradient state in subregion in $\omega$ of $\Omega$ in general case only.

### 3. Sensor and Asymptotic $\omega_G$-Detectability:

This concept of regional detectability is introduced and extended by Al-Saphory and El-Jai as in (Al-Saphory and El Jai, 2002), (Al-Saphory, 2011) and (Al-Saphory *et al.*, 2010) and then for the boundary case see (Al-Saphory, 2002) and (Al-Saphory *et al.*, 2016). In this section, we extend some definitions and characterizations to the concept of asymptotic regional gradient detectability in order to construct an general asymptotic $\omega_G$-observer for the current gradient state of the original system. For this purpose we need some definitions and characterization as in (Al-Saphory *et al.*, 2013) and (Al-Saphory *et al.*, 2015).

### 3.1 Definitions and Characterizations:
***Definition 3.1:***
The systems (1)-(2) are said to be exactly regionally gradient observable on $\omega$ (exactly $\omega_G$ –observable) if

$Im\ H = Im\ \chi_\omega \nabla K^* = (H^1(\omega))^n$

*Definition 3.2:*

The systems (1)–(2) are said to be weakly regionally gradient observable on $\omega$ (weakly $\omega_G$ –observable) if $\overline{ImH} = \overline{Im\ \chi_\omega \nabla K^*} = (H^1(\omega))^n$

*Remark 3.3:*

The definition 3.2 is equivalent to say that the systems (1)-(2) are weakly $\omega_G$ –observable if $ker\ H^* = ker\ K\nabla^* \chi_\omega^* = \{0\}$

*Definition 3.4:*

A sensor $(D, f)$ is gradient strategic on $\omega$ ($\omega_G$-strategic) if the observed system is weakly $\omega_G$-observable.

*Definition 3.5:*

The system (1) is said to be gradient stable ($\Omega_G$-stable) if the operator $A$ generates a semi-group which is gradient stable on the $(H^1(\Omega))^n$. It is easy to see that the system (1) is $\Omega_G$-stable, if and only if for some positive constants $M, \alpha$, we have
$\|\nabla S_A(.)\|_{(H^1(\Omega))^n} \le Me^{-\alpha t}, \forall t \ge 0$

If $(S_A(t))_{t\ge 0}$ is $\Omega_G$-stable semi-group in $(H^1(\Omega))^n$, then for all $x_0 \in H^1(\Omega)$, the solution of associated system satisfies
$\lim_{t\to\infty}\|\nabla x(.,t)\|_{(H^1(\Omega))^n} = \lim_{t\to\infty}\|\nabla S_A(.)x_0\|_{(H^1(\Omega))^n} = 0$ (5)

*Definition 3.6:*

The systems (1)-(2) are said to be asymptotic gradient detectable (asymptotic $\Omega_G$-detectable) if there exists an operator $H_\Omega: R^q \to (H^1(\Omega))^n$ such that $(A - H_\Omega C)$ generates a strongly continuous semi-group $(S_{H_\Omega}(t))_{t\ge 0}$ which is asymptotically $G$-stable on $(H^1(\Omega))^n$.

*Definition 3.7:*

The system (1) is said to be regional gradient stable ($\omega_G$-stable) if the operator $A$ generates a semi-group which is gradient stable on the $(H^1(\omega))^n$. It is easy to see that the system (1) is $\omega_G$-stable, if and only if for some positive constants $M_{\omega_G}, \alpha_{\omega_G}$, we have
$\|\chi_\omega \nabla S_A(.)\|_{(H^1(\omega))^n} \le M_{\omega_G} e^{-\alpha_{\omega_G} t}, \forall t \ge 0$

If $(S_A(t))_{t\ge 0}$ is $\omega_G$-stable semi-group in $(H^1(\omega))^n$, then for all $x_0 \in H^1(\Omega)$, the solution of associated system satisfies
$\lim_{t\to\infty}\|\chi_\omega \nabla x(.,t)\|_{(H^1(\omega))^n} = \lim_{t\to\infty}\|\chi_\omega \nabla S_A(.)x_0\|_{(H^1(\omega))^n} = 0$ (6)

*Definition 3.8:*

The systems (1)-(2) are said to be asymptotic regional gradient detectable (asymptotic $\omega_G$-detectable) if there exists an operator $H_{\omega_G}: R^q \to (H^1(\omega))^n$ such that $(A - H_{\omega_G} C)$ generates a strongly continuous semi-group $(S_{H_{\omega_G}}(t))_{t\ge 0}$ which is asymptotically $G$-stable on $(H^1(\omega))^n$.

*Remark 3.9:*

In this paper, we only need the relation (6) to be true on a sub-region $\omega$ of the region $\Omega$
$\lim_{t\to\infty}\|\nabla x(.,t)\|_{(H^1(\omega))^n} = 0$

Now, to study the relation between the asymptotic $\omega_G$-detectable and $\omega_G$-strategic sensors. Let us consider the systems (1)-(2) which are observed by $q$ sensors $(D_i, f_i)_{1\le i\le q}$ with $D_i \subset \Omega$ and $f_i \in H^1(\Omega)$ for $i = 1, ..., q$. We assume that the operator $A$ has a complete set of eigenfunctions denoted by $\psi_{mj}$ in $(H^1(\Omega))^n$ orthonormal in $(H^1(\omega))^n$ associated with the eigenvalues $\lambda_m$ of multiplicity $r_m$ and suppose that the system (1) has unstable modes. We have the following sufficient condition for existence of asymptotic $\omega_G$-detectability in terms of the strategic sensors structure.

*Theorem 3.10:*

Suppose that there are $q$ zone sensors $(D_i, f_i)_{1\le i\le q}$ and the spectrum of $A$ contains $J$ eigenvalues with non-negative real parts. The system (1) together with output function (2) is $\omega_G$- detectable if and only if :
1. $q \ge r$
2. rank $G_m = r_m$ , $\forall m, m = 1, ..., J$ with

$$G_m = (G_m)_{ij} = \begin{cases} \psi_{m_j}(b_i), f_i(.) >_{L^2(D_i)} \text{ for zone sensors} \\ \psi_{m_j}(b_i) \text{ for pointwise sensors} \\ < \frac{\partial \psi_{m_j}}{\partial v}, f_i(.) >_{L^2(\Gamma_i)} \text{ for boundary zone sensors} \end{cases}$$

where $\sup r_m = m < \infty$ and $j = 1, \ldots, r_m$.

**Proof:**

The proof is limited to the case of zone sensors. Under the assumptions of section 2, the system (1) can be decomposed by the projections $P$ and $I - P$ on two parts, unstable and stable. The state vector may be given by $x(\xi, t) = [x_1(\xi, t), x_2(\xi, t)]^{tr}$ where $x_1(\xi, t)$ is the state component of the unstable part of the system (1), may be written in the form

$$\begin{cases} \frac{\partial x_1}{\partial t}(\xi, t) = A_1 x_1(\xi, t) + PBu(t) & Q \\ x_1(\xi, 0) = x_{01}(\xi) & \Omega \\ x_1(\eta, t) = 0 & \Theta \end{cases} \quad (7)$$

and $x_2(\xi, t)$ is the component state of the stable part of the system (1) given by

$$\begin{cases} \frac{\partial x_2}{\partial t}(\xi, t) = A_2 x_2(\xi, t) + (I - P)Bu(t) & Q \\ x_2(\xi, 0) = x_{02}(\xi) & \Omega \\ x_2(\eta, t) = 0 & \Theta \end{cases} \quad (8)$$

The operator $A_1$ is represented by a matrix of order $(\sum_{m=1}^{J} r_m, \sum_{m=1}^{J} r_m)$ given by $A_1 = \text{diag}[\lambda_1, \ldots, \lambda_1, \lambda_2, \ldots, \lambda_2, \ldots, \lambda_J, \ldots, \lambda_J]$ and $PB = [G_1^{tr}, G_2^{tr}, \ldots, G_J^{tr}]$. The condition (2) of this theorem, allows that the suit $(D_i, f_i)_{1 \leq i \leq q}$ of sensors is $\omega_G$-strategic for the unstable part of the system (1), the subsystem (7) is weakly $\omega_G$-observable(Al-Saphory *et al.*, 2015), and since it is finite dimensional, then it is exactly $\omega_G$-observable. Therefor it is $\omega_G$-detectable, and hence there exists an operator $H_\omega^1$ such that $(A_1 - H_\omega^1 C)$ which is satisfied the following:

$\exists M_\omega^1, \alpha_\omega^1 > 0$ such that $\|e^{(A_1 - H_\omega^1 C)t}\|_{(H^1(\omega))^n} \leq M_\omega^1 e^{-\alpha_\omega^1(t)}$

and then we have

$\|x_1(., t)\|_{(H^1(\omega))^n} \leq M_\omega^1 e^{-\alpha_\omega^1(t)} \|Px_0(.)\|_{(H^1(\omega))^n}$.

Since the semi-group generated by the operator $A_2$ is stable on $(H^1(\omega))^n$, then there exist $M_\omega^2, \alpha_\omega^2 > 0$ such that

$\|x_2(., t)\|_{(H^1(\omega))^n} \leq M_\omega^2 e^{-\alpha_\omega^2(t)} \|(I - P)x_0(.)\|_{(H^1(\omega))^n}$

$+ \int_0^t M_\omega^2 e^{-\alpha_\omega^2(t)} \|(I - P)x_0(.)\|_{(H^1(\omega))^n} \|u(\mathcal{T})\| d\mathcal{T}$

and therefor $x(\xi, t) \to 0$ when $t \to \infty$. Finally, the system (1)-(2) are $\omega_G$-detectable (Al-Saphory *et al.*, 2010).

Reciprocally, if the system (1) together with the output function (2) is $\omega_G$-detectable, there exists an operator $H_\omega \in L(R^q, (H^1(\omega))^n)$, such that $(A - H_{\omega_G} C)$ generates a $\omega_G$-stable, strongly continuous semi-group $(S_{H_\omega}(t))_{t \geq 0}$ on the space $(H^1(\omega))^n$, which is satisfied the following:

$\exists M_\omega, \alpha_\omega > 0$ such that $\|\chi_\omega \nabla S_{H_\omega}(.)\|_{(H^1(\omega))^n} \leq M_\omega e^{-\alpha_\omega(t)}$.

Thus, the unstable subsystem (7) is $\omega_G$-detectable. We recall that a system is weakly $\omega_G$-observable, i.e. $[K \nabla^* \chi_\omega^* x^*(., t) = 0 \Rightarrow x^*(., t) = 0]$ (Al-Saphory *et al.*, 2015). For $x^*(., t) \in (H^1(\omega))^n$, we have

$$K \nabla^* \chi_\omega^* x^*(., t) = (\sum_{m=1}^{J} e^{\lambda_{m_j} t} < \psi_{m_j}(.), \chi_\omega^* \nabla x^*(., t) >_{(H^1(\omega))^n}, < \psi_{m_j}(.), f_j(., t) >_{H^1(\Omega)})_{1 \leq i \leq q}$$

$$= \sum_{m=1}^{J} e^{\lambda_{m_j} t} < \nabla \psi_{m_j}(.), x^*(., t) >_{(H^1(\omega))^n} < \psi_{m_j}(.), f_j(., t) >_{H^1(\Omega)})_{1 \leq i \leq q}$$

$$= \sum_{m=1}^{J} e^{\lambda_{m_j} t} < \psi_{m_j}(.), x^*(., t) >_{(H^1(\omega))^n} < \psi_{m_j}(.), f_j(., t) >_{(H^1(\omega))^n})_{1 \leq i \leq q}$$

If the rank $G_m x_m \neq r_m$ for $m, m = 1, \ldots, J$, there exists $x^*(., t) \in (H^1(\omega))^n$, such that $K \nabla^* \chi_\omega^* x^*(., t) = 0$, this leads

$\sum_{m=1}^{J} < \psi_{m_j}(.), x^*(., t) >_{(H^1(\omega))^n} < \psi_j(.), f_j(., t) >_{H^1(\Omega)} = 0$.

The state vectors $x_m$ may be given by

$x_m(., t) = [< \psi_{1_j}(.), x^*(., t) >_{(H^1(\omega))^n} < \psi_{J_m}(.), x^*(., t) >_{H^1(\Omega)}]^{tr} \neq 0$

we then obtain $G_m x_m = 0$ for $m, m = 1, ..., J$. Consequently, the subsystem (7) is not weakly $\omega_G$-observable and therefore the suite $(D_i, f_i)_{1 \leq i \leq q}$ of sensors is not $\omega_G$-strategic. Thus, the system (1)-(2) is not $\omega_G$-detectable. Finally we have the rank $G_m \neq r_m$ for all $m, m = 1, ..., J$.□

Now, it is clear that:
1. A system which is $G$- detectable, is $\omega_G$- detectable.
2. A system which is exponentially $\omega_G$- detectable, is asymptotically $\omega_G$- detectable.
3. A system which is asymptotically $\omega_G$-detectable, is $\omega_G^1$- detectable, for every subset $\omega_1$ of $\omega$ but the converse is not true.

## 4. General Asymptotic $\omega_G$-observers:

The purpose of this section is to introduced an approach which enable to achieve the existence of general asymptotic regional gradient observer ($GA\omega_G$-observer) in general case is derived from (El-Jai and Amouroux, 1988), (Al-Saphory and El-Jai, 2002) and (Al-Saphory *et al*., 2010).This approach which allows to construct the current gradient state in $\omega$ of the systems (1)-(2) in connection of asymptotically $\omega_G$-detectability and $\omega_G$-strategic sensors in order to characterize the asymptotic regional gradient observability. Thus, it's necessary to approximate the dynamic characteristics of the observer in subspace spanned by a finite number of suitably chosen bases. Let us consider the operator $A$ has a complete set of eigenfunctions $\varphi_{nj}$ in $(H^1(\Omega))^n$ orthonormal to $(H^1(\omega))^n$ associated with the eigenvalues $\lambda_n$ of multiplicity $r_n$ and suppose that the system (1) has unstable modes. Then, we present sufficient conditions for existing a $GA\omega_G$-observer in order to provide an approximation to the gradient state of observed system.

### 4.1 Definitions:
**Definition 4.1:**

Suppose there exists a dynamical system with state $z(.,t) \in Z$ given by

$$\begin{cases} \frac{\partial z(.,t)}{\partial t} = F_{\omega_G} z(.,t) + G_{\omega_G} u(t) + H_{\omega_G} y(t) & Q \\ z(\xi, 0) = z_0(\xi) & \Omega \\ z(\eta, t) = 0 & \Sigma \end{cases} \quad (9)$$

Where $F_{\omega_G}$ generator a strongly continuous semi-group $(S_{F_{\omega_G}}(t))_{t \geq 0}$ on separable Hilbert space $Z$ which is asymptotically $\omega_G$ –stable. Thus, $\exists M_{F_{\omega_G}}, \alpha_{F_{\omega_G}} > 0$ such that

$$\left\| S_{F_{\omega_G}}(.) \right\| \leq M_{F_{\omega_G}} e^{-\alpha_{F_{\omega_G}} t}, \forall t \geq 0.$$

and let $G_{\omega_G} \in L(U, Z), H_{\omega_G} \in L(\mathcal{O}, Z)$ such that the solution of (9) similar to (3)

$$z(\xi, t) = S_{F_{\omega_G}}(t) z(\xi) + \int_0^t S_{F_{\omega_G}}(t - \tau) [G_{\omega_G} u(\tau) + H_{\omega_G} y(\tau)] d\tau$$

Thesystem (9) defines asymptotic regional gradient estimator for $\chi_\omega \nabla T x(\xi, t)$ where $x(\xi, t)$ is the solution of the systems (1)-(2) if

$\lim_{t \to \infty} \| z(.,t) - \chi_\omega \nabla T x(\xi, t) \|_{(H^1(\omega))^n} = 0$

and $\chi_\omega \nabla T$ maps $D(A)$ into $D(F)$ where $z(\xi, t)$ is the solution of system (9). And (9) specifies an asymptotic $\omega_G$-observer of the system given by (1) and (2) if the following holds:
1- There exists $M_{\omega_G} \in L(R^q, (H^1(\omega))^n)$ and $N_{\omega_G} \in L((H^1(\omega))^n)$ such that $M_{\omega_G} C + N_{\omega_G} \chi_\omega \nabla T = I_{\omega_G}$.
2- $\chi_\omega \nabla T A - F_{\omega_G} \chi_\omega \nabla T = H_{\omega_G} C$, $G_{\omega_G} = \chi_\omega \nabla T B$.
3- The system (9) defines an asymptotic $\omega_G$-estimator for $\chi_\omega \nabla T x(\xi, t)$.

The object of an asymptotic $\omega_G$-estimator (asymptotic $\omega_G$-observer) is to provide an approximation to the gradient state of the original system. This approximation is given by
$\hat{x}(t) = M_{\omega_G} y(t) + N_{\omega_G} z(t)$

**Definition 4.2:**

The systems (1)-(2) are general asymptotic regional gradient observable ($AG\omega_G$-observable) if there exists a dynamical system which is $GA\omega_G$-observer for the original system.

### 4.2 $GA\omega_G$-observer reconstruction method:

The problem of studying asymptotic $\omega_G$-observability may be through the observation operator $C$, that means, we can the characterize the $GA\omega_G$-observer by a good choice of the sensors structure. For this objective suppose the system has $q$ sensors $(D_i, f_i)_{1 \leq i \leq q}$.

Consider the system

$$\begin{cases} \frac{\partial x}{\partial t}(\xi,t) = Ax(\xi,t) + Bu(t) & Q \\ x(\xi,0) = x_0(\xi) & \Omega \\ x(\eta,t) = 0 & \Sigma \end{cases} \qquad (10)$$

In this case the output function (2)

$$y(t) = Cx(.,t) \qquad (11)$$

Let $\omega$ be a given subdomain of $\Omega$ and suppose that $T \in \mathcal{L}((H^1(\Omega))^n)$, and $\chi_\omega \nabla Tx(\xi,t) = T_\omega \nabla x(\xi,t)$ there exists a system with state $z(\xi,t)$ such that

$$z(\xi,t) = \chi_\omega \nabla Tx(\xi,t) \qquad (12)$$

From equation (11) and (12) we have

$$\begin{bmatrix} y \\ z \end{bmatrix} = \begin{bmatrix} C \\ \chi_\omega \nabla T \end{bmatrix} x$$

If we assume that there exist two bounded linear operators $M_{\omega_G}: \mathcal{O} \to ((H^1(\omega))^n$ and $N_{\omega_G}: ((H^1(\omega))^n \to ((H^1(\omega))^n$, such that $M_{\omega_G}C + N_{\omega_G}T_\omega = I$, then by deriving $z(\xi,t)$ in (12) we have

$$\frac{\partial z}{\partial t}(\xi,t) = \chi_\omega \nabla T \frac{\partial x}{\partial t}(\xi,t) = \chi_\omega \nabla TAx(\xi,t) + \chi_\omega \nabla TBu(t)$$
$$= \chi_\omega \nabla TAM_{\omega_G}y(\xi,t) + \chi_\omega \nabla TAN_{\omega_G}z(\xi,t) + \chi_\omega \nabla TBu(t)$$

Now, consider the system (which is destined to be the $\omega_G$-observer)

$$\begin{cases} \frac{\partial \hat{z}}{\partial t}(\xi,t) = F_{\omega_G}\hat{z}(\xi,t) + G_{\omega_G}u(t) + H_{\omega_G}y(.,t) & Q \\ \hat{z}(\xi,0) = z_0(\xi) & \Omega \\ \hat{z}(\eta,t) = 0 & \Sigma \end{cases} \qquad (13)$$

Where $F_{\omega_G}$ generates a strongly continuous semi-group $(S_{F_{\omega_G}}(t))_{t\geq 0}$ on separable Hilbert space $Z$ which is asymptotically gradient stable thus,

$\exists M_{F_{\omega_G}}, \alpha_{F_{\omega_G}} > 0$ such that $\left\| \chi_\omega \nabla S_{F_{\omega_G}}(.) \right\| \leq M_{F_{\omega_G}} e^{-\alpha_{F_{\omega_G}}t}, \forall t \geq 0$.

and let $G_{\omega_G} \in L(U,Z)$, $H_{\omega_G} \in L(\mathcal{O},Z)$ such that the solution of (13) is given by

$$\hat{z}(\xi,t) = S_{F_{\omega_G}}(t)\hat{z}_0(\xi) + \int_0^t S_{F_{\omega_G}}(t-\tau)[G_{\omega_G}u(\tau) + H_{\omega_G}x(b_i,\tau)]d\tau$$

Now, we present the main result which enable to observe asymptotically the current gradient state in $\omega$, by the gradient state of the system (13).

***Theorem 4.3:***

Suppose that the operator $F_{\omega_G}$ generates a strongly continuous semi-group which is $\omega_G$-stable on $((H^1(\omega))^n$, then, the system (13) is $GA\omega_G$-observer for systems (1)-(2), that is,

$\lim_{t\to\infty}[\chi_\omega \nabla Tx(\xi,t) - \hat{z}(\xi,t)] = 0$

If the following conditions hold:

1- There exists $M_{\omega_G} \in L(R^q, (H^1(\omega))^n)$ and $N_{\omega_G} \in L((H^1(\omega))^n)$ such that $M_{\omega_G}C + N_{\omega_G}\chi_\omega \nabla T = I_\omega$.

2- $\chi_\omega \nabla TA - F_{\omega_G}\chi_\omega \nabla T = H_{\omega_G}C$, $G_{\omega_G} = \chi_\omega \nabla TB$.

***Proof:***

Let $z(\xi,t) = \chi_\omega \nabla Tx(\xi,t)$ and $\hat{z}(\xi,t)$ be a solution of (13) and let us denote the observer error by the following form

$e(\xi,t) = z(\xi,t) - \hat{z}(\xi,t)$

We have

$$\frac{\partial e}{\partial t}(\xi,t) = \frac{\partial z}{\partial t}(\xi,t) - \frac{\partial \hat{z}}{\partial t}(\xi,t) = \chi_\omega \nabla TAx(\xi,t) + \chi_\omega \nabla TBu(t) - F_{\omega_G}\hat{z}(\xi,t) - G_{\omega_G}u(t) - H_{\omega_G}y(\xi,t)$$
$$= F_{\omega_G}e(\xi,t) - F_{\omega_G}z(\xi,t) + \chi_\omega \nabla TAx(\xi,t) - H_{\omega_G}y(\xi,t) + \chi_\omega \nabla TBu(t) - G_{\omega_G}u(t)$$
$$= F_{\omega_G}e(\xi,t) + [\chi_\omega \nabla TAx(\xi,t) - F_{\omega_G}\chi_\omega \nabla Tx(\xi,t) - H_{\omega_G}Cx(\xi,t)] + [\chi_\omega \nabla TBu(t) - G_{\omega_G}u(t)]$$
$$= F_{\omega_G}e(\xi,t) + [\chi_\omega \nabla TA - F_{\omega_G}\chi_\omega \nabla T - H_{\omega_G}C]x(\xi,t) + [\chi_\omega \nabla TB - G_{\omega_G}]u(t) = F_{\omega_G}e(\xi,t)$$

Consequently $e(\xi,t) = S_{F_{\omega_G}}(t)[\chi_\omega \nabla Tx_0(.) - \hat{z}_0(.)]$. the $G$-stability of the operator allows to obtain

$\|e(\xi,t)\|_{(H^1(\omega))^n)} \leq M_{F_{\omega_G}}e^{-\alpha_{F_{\omega_G}}t}\|\chi_\omega \nabla Tx_0(\xi) - \hat{z}_0(\xi)\|_{(H^1(\omega))^n)}$

And therefore $\lim_{t\to\infty} e(\xi,t) = 0$.

Now, let the approximate solution to the gradient state of the original system is $\hat{x}(\xi,t) = M_{\omega_G} y(.,t) + N_{\omega_G} \hat{z}(\xi,t)$, then we have

$$\hat{e}(\xi,t) = x(\xi,t) - \hat{x}(\xi,t) = x(\xi,t) - M_{\omega_G} y(.,t) - N_{\omega_G} \hat{z}(\xi,t),$$
$$= x(\xi,t) - M_{\omega_G} C x(\xi,t) - N_{\omega_G} \chi_\omega \nabla T x(\xi,t) + N_{\omega_G} [\chi_\omega \nabla T x(\xi,t) - \hat{z}(\xi,t)]$$
$$= N_{\omega_G} [\chi_\omega \nabla T x(\xi,t) - \hat{z}(\xi,t)] = N_{\omega_G} [z(\xi,t) - \hat{z}(\xi,t)] = N_{\omega_G} e(\xi,t)$$

Finally, we have
$\lim_{t\to\infty} \hat{z}(\xi,t) = z(\xi,t)$. □

Thus, we can deduced a sufficient condition for existence $GA\omega_G$-observer is formulated in the following theorem.

*Remark 4.4:*

We can deduced that
1. The theorem 4.3 gives the conditions which guarantee that the dynamical system (9) is a $GA\omega_G$-observer for the system (1)-(2).
2. From theorem 4.3 we get the relation between the regional gradient strategic sensor and $GA\omega_G$-observer.
3. A system which is an $AGG$-observer is $AG\omega_G$-observer.
4. If a system is $GA\omega_G$-observer, then it is $GA\omega_G^1$-observer in every subset $\omega_1$ of $\omega$, but the converse is not true. This is may be explain in the following example;

*Example 4.5:*

Consider the system

$$\begin{cases} \frac{\partial x}{\partial t}(\xi,t) = \gamma_1 \frac{\partial^2 x}{\partial \xi^2}(\xi,t) + \gamma_2 x(\xi,t) & ]0,a[, t>0 \\ x(\xi,0) = x_0(\xi) & ]0,a[ \\ x(0,t) = x(a,t) = 0 & t>0 \end{cases} \quad (14)$$

Where $\gamma_1 > 0$, $\gamma_2 > 0$ and $\Omega = ]0,a[$. Let $b_i \in \Omega$ are the locations of the pointwise sensors $(b_i, \delta_{b_i})$. Then the augmented output function is given by:

$$\mathbf{y}(.,\mathbf{t}) = \int_\Omega x(\xi,t)\delta(\xi - b_i)d\xi. \quad (15)$$

Now, consider the dynamical system

$$\begin{cases} \frac{\partial z}{\partial t}(\xi,t) = \gamma_1 \frac{\partial^2 z}{\partial \xi^2}(\xi,t) + \gamma_2 z(\xi,t) - HC(z(\xi,t) - x(\xi,t)) & ]0,a[, t>0 \\ z(\xi,0) = z_0(\xi) & ]0,a[ \\ z(\eta,t) = 0 & t>0 \end{cases} \quad (16)$$

Let $\omega = (\alpha, \beta)$ be a subregion of $\Omega$. The eigenfunctions and the eigenvalue related to the operator $(\gamma_1 \frac{\partial^2}{\partial \xi^2} + \gamma_2)$ are given by:

$\varphi_n(\xi) = \left(\frac{2}{\beta-\alpha}\right)^{\frac{1}{2}} sinn\pi(\frac{\xi-\alpha}{\beta-\alpha})$ and $\lambda_n = 1 - \left(\frac{n\pi}{\beta-\alpha}\right)^2$.

Now, if $\frac{b_i}{a} \in Q \cap (0,a)$ then the sensor is not $\omega_G$-strategic (Al-Saphory *et al*., 2015) for unstable subsystem of (14) and therefor the system (14)-(15) is not $\omega_G$-detectable in $\Omega$ (Al-Saphory *et al*., 2010). Then the dynamical system (16) is not $\omega_G$ – observer for the system (14)-(15). Thus, the system (14)-(15) is $\omega_G$-detectable if $(\frac{b_i - \alpha}{\beta - \alpha}) \notin Q \cap (0,a)$. That means, the sensors $(b_i, \delta_{b_i})$ are $\omega_G$-strategic to unstable subsystem of (14) and therefor the systems (14)-(15) are $\omega_G$-detectable. Therefor, the dynamical systems (16) are $\omega_G$ – observer for the system (14)-(15) (Al-Saphory *et al*., 2010). □

*5. Application to $GA\omega_G$-observer:*

In this section, we give some results related to different types of measurements domains and boundary conditions. We consider the distributed diffusion systems defined on two dimensional domain $]0,a_1[\times]0,a_2[$ on $\Omega$.

$$\begin{cases} \frac{\partial x}{\partial t}(\xi_1,\xi_2,t) = \frac{\partial^2 x}{\partial \xi_1^2}(\xi_1,\xi_2,t) + \frac{\partial^2 x}{\partial \xi_2^2}(\xi_1,\xi_2,t) & Q \\ x(\xi_1,\xi_2,0) = x_0(\xi_1,\xi_2) & \Omega \\ x(\xi,\eta,t) = 0 & \Sigma \end{cases} \quad (17)$$

The augmented output function is given by $y(.,t) = \int_{\Gamma_0} x(\eta_1, \eta_2, t) f(\eta_1, \eta_2) d\eta_1 d\eta_2$ (18) Let $\Gamma = \{a_1\} \times ]0, a_2[$ be a region on $]0, a_1[ \times ]0, a_2[$. The eigenfunctions of the operator $(\frac{\partial^2}{\partial \xi_1^2} + \frac{\partial^2}{\partial \xi_2^2})$ for the Dirichlet boundary condition are defined by

$$\varphi_{ij}(\xi_1, \xi_2) = \frac{2}{\sqrt{a_1 a_2}} \sin i\pi \frac{\xi_1}{a_1} \sin j\pi \frac{\xi_2}{a_2}$$

Associated with the eigenvalues

$$\lambda_{ij} = (\frac{i^2}{a_1^2} + \frac{j^2}{a_2^2})\pi^2$$

(The dynamical system

$$\begin{cases} \frac{\partial z}{\partial t}(\xi_1, \xi_2, t) = \frac{\partial^2 z}{\partial \xi_1^2}(\xi_1, \xi_2, t) + \frac{\partial^2 z}{\partial \xi_2^2}(\xi_1, \xi_2, t) + Bu(\xi_1, \xi_2, t) - H_{\omega_G}(C z(\xi_1, \xi_2, t) - y(t) & Q \\ z(\xi_1, \xi_2, 0) = z_0(\xi_1, \xi_2) & \Omega \\ x(\eta_1, \eta_2, t) = 0 & \Sigma \end{cases}$$ (19)

together with the system (17)– (18) are equivalent to the systems (1)–(2) and (9) with the operator $A = \left(\frac{\partial^2}{\partial \xi_1^2} + \frac{\partial^2}{\partial \xi_2^2}\right)$.

***Zone sensor cases:***

Let the measurement support is rectangular with $]\xi_1 - l_1, \xi_1 + l_1[ \times ]\xi_2 - l_2, \xi_2 + l_2[ \in \Omega$ (see Figure 3).

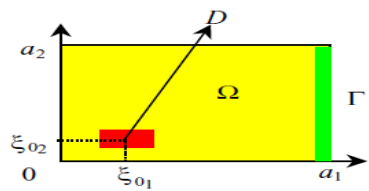

**Fig. 3:** Location of internal zone sensor $D$.

If $f_1$ is symmetric about $\xi_1 = \xi_{01}$ and $f_2$ is symmetric with respect to $\xi_2 = \xi_{02}$ then we have get the following result:

***Corollary 5.1:***

The dynamical system (19) is $AG\omega_G$-observer for the system (17)-(18) if $\frac{i\xi_{01}}{a_1}$ and $\frac{j\xi_{02}}{a_2} \notin N, \forall i, j = 1, \dots, J$.

In the case where $\Gamma \subset \partial\Omega$ and $f \in L^2(\Gamma)$, the sensor $(D, f)$ may be located on the boundary in $\Gamma_0 = ]\eta_{01} - l_1, \eta_{01} + l_1[ \times \{a_2\}$, then we have:

***Corollary 5.2:***
***One side case (see Figure 4):***

suppose that the sensor $(D, f)$ is located on $\Gamma_0 = ]\eta_{01} - l_1, \eta_{01} + l_1[ \times \{a_2\} \subset \partial\Omega$, and $f$ is symmetric with respect to $\eta_1 = \eta_{01}$ then the dynamical systems (19) are $GA\omega_G$-observer for the system (17)-(18) if $\frac{i\xi_{01}}{a_1} \notin N, \forall i, j = 1, \dots, J$.

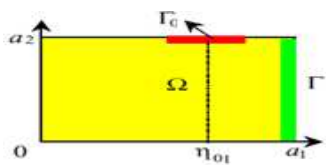

**Fig. 4:** Location of boundary zone sensor one side case

*Two side case (see Figure 5):*
    suppose that the sensor $(D, f)$ is located on $\varGamma_0 = ]0, \eta_{01} + l_1[ \times \{0\} \cup \{0\} \times ]0, \eta_{02} + l_2[ \subset \partial\Omega$, and $f|_{\varGamma_1}$ is symmetric with respect to $\eta_1 = \eta_{01}$ and the function $f|_{\varGamma_2}$ is symmetric with respect to $\eta_2 = \eta_{02}$, then the dynamical system (19) is $GA\omega_G$-observerfor the system (17)-(18) if $\frac{i\xi_{01}}{a_1}$ and $\frac{j\xi_{02}}{a_2} \notin N, \forall i,j = 1, \dots, J$.

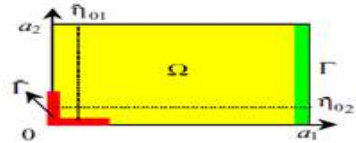

**Fig. 5:** Location of boundary zone sensor two side case

*ii. Pointwise sensor cases:*
    In this cases if if $b = (b_1, b_2) \in \partial\Omega$ then, we have:

*Corollary 5.3:*
*Internal case (see Figure 6):*
    The dynamical system(19) is $GA\omega_G$-observer for the systems(17)-(18) if $\frac{ib_1}{a_1}$ and $\frac{jb_2}{a_2} \notin N, \forall i,j = 1, \dots, J$.

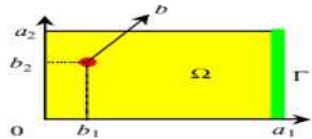

**Fig. 6:** Location of Internal pointwise sensors.

*Filament case (see Figure 7):*
    Suppose that the observation is given by the filament sensor where $\sigma = Im(\gamma)$ is symmetric with respect to the line $b = (b_1, b_2)$, if $\frac{ib_1}{a_1}$ and $\frac{jb_2}{a_2} \notin N, \forall i,j = 1, \dots, J$. Then the dynamical system (19) is $GA\omega_G$-observer for the systems(17)-(18).

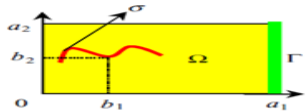

**Fig. 7:** Location of Internal pointwise sensors (filament case).

*Boundary case (see Figure 8):*
    The dynamical system (19) is $GA\omega_G$-observerfor the systems(17)-(18) if $\frac{jb_2}{a_2} \notin N, \forall j = 1, \dots, J$.

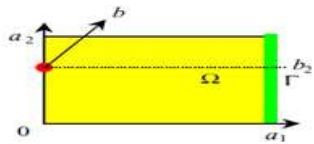

**Fig. 8:** Location of boundary pointwise sensor b.

*Remark 5.4:*
    We can extend these results to the case where $\omega$ is sub-region of the boundary of the domain $\Omega$ as in(Al-Saphory *et al*., 2016).

*Conclusion:*

The concept have been studied in this paper is related to the existing and characterizing an $GA\omega_G$-observer only for a distributed parameter systems. More precisely, we have given an approach for building a $GA\omega_G$-estimatorwhich reconstruct a gradient state in considered sub-region $\omega$. Also, we show that there exists a dynamical system is not $GA\omega_G$-observer in usual sense but it is $GAG$-observer. For the future work, one can extension of these results to the problem in identity case and in reduced order case as in (Al-Saphory and Al-Mullah, 2015Al-Saphory and Jaafar 2015).

## ACKNOWLEDGMENTS

Our thanks in advance to the editors and experts for considering this paper to publish in this estimated journal and for their efforts.